\documentclass[doublecolumn]{IEEEtran}
\usepackage{verbatim,amsmath,amsfonts,amssymb,mathrsfs,color,graphicx,algorithm,algorithmic,cite,subfigure}
\newtheorem{theorem}{Theorem}
\newtheorem{lemma}{Lemma}
\newtheorem{proposition}{Proposition}

\newtheorem{assumption}{Assumption}

\newenvironment{remark}[1][Remark]{\begin{trivlist}
\item[\hskip \labelsep {\bfseries #1}]}{\end{trivlist}}

\newcommand{\qed}{\nobreak \ifvmode \relax \else
      \ifdim\lastskip<1.5em \hskip-\lastskip
      \hskip1.5em plus0em minus0.5em \fi \nobreak
      \vrule height0.75em width0.5em depth0.25em\fi}
\DeclareMathOperator*{\argmin}{arg\,min}
\DeclareMathOperator*{\argmax}{arg\,max}
\def\dist{\mathrm{dist}}

\def\proj{\mathsf{\Pi}}

\def\a{\alpha}

\def\g{\gamma}

\def\o{\omega}
\def\O{\Omega}

\def\Gc{\mathcal{G}}
\def\Vc{\mathcal{V}}
\def\Ec{\mathcal{E}}

\def\Xc{\mathcal{X}}
\def\Fc{\mathcal{F}}
\def\Nc{\mathcal{N}}
\def\Es{\mathsf{E}}

\def\1b{\mathbf{1}}
\def\0b{\mathbf{0}}

\def\la{\langle}
\def\ra{\rangle}

\def\diag{\textrm{diag}}

\def\Rbb{\mathbb{R}}

\def\sml#1{{#1}}
\def\smlr#1{{#1}}

\allowdisplaybreaks
\newcommand*{\medcap}{\mathbin{\scalebox{1.5}{\ensuremath{\cap}}}}
\title{Approximate Projection Methods for Decentralized Optimization with Functional Constraints}
\author{Soomin Lee and Michael M. Zavlanos\thanks{
This work is supported by ONR under grant \#N000141410479.
The research was conducted when Soomin Lee was a Postdoctoral Associate with the Dept. of Mechanical Engineering and Materials Science, Duke University,
Durham, NC, 27708, USA.
Michael M. Zavlanos is with the Dept. of
Mechanical Engineering and Materials Science, Duke University, Durham,
NC, 27708, USA. e-mails: \texttt{soominl@gmail.com, michael.zavlanos@duke.edu.}
}
}

\begin{document}

\maketitle

\begin{abstract}
We consider distributed convex optimization problems that involve a separable objective function and nontrivial functional constraints, such as Linear Matrix Inequalities (LMIs). We propose a decentralized and computationally inexpensive algorithm %to solve such problems over time-varying directed networks of agents,
which is based on the concept of approximate projections.
Our algorithm is one of the consensus based methods in that, at every iteration, each agent performs a consensus update of its decision variables followed by an optimization step of its local objective function and
local constraints.
Unlike other methods, the last step of our method is not an Euclidean projection onto the feasible set, but instead a subgradient step in the direction that minimizes the local constraint violation.
%This way, our method avoids expensive computations, such as eigenvalue decompositions when projecting onto semidefinite constraints. %, significantly reducing computational cost.
We propose two different averaging schemes to mitigate the disagreements among the agents' local estimates.
We show that the algorithms converge almost surely, i.e., every agent agrees on the same optimal solution, under the assumption that the objective functions and constraint functions are nondifferentiable and their subgradients are bounded.
%We also establish the finite termination of the algorithm in the case of a distributed set feasibility problem.
We provide simulation results on a decentralized optimal gossip averaging problem, which involves SDP constraints, to complement our theoretical results.
%and a decentralized feasibility (robust LQR) problem, both involving large numbers of LMI constraints.
\end{abstract}

\section{Introduction}
\IEEEPARstart{D}{ecentralized} optimization has been extensively studied in recent years due to
a variety of applications in machine learning, signal processing, and control
for robotic networks, sensor networks, power networks, and wireless communication networks %\cite{Langford-book,Bianchi-book, Camponogara, johansson,ram_info,rabbat,con03,con01}.
\cite{Bianchi-book, Camponogara,ram_info,con03,con01}.
A number of problems arising in these areas can be cast as distributed convex optimization problems over multiagent networks,
where individual agents cooperatively try to minimize a common cost function over a common constraint set
in the absence of full knowledge about the global problem structure.
The main feature of carrying these optimizations over networks
is that the agents can only communicate with their neighboring agents.
This communication structure can be cast as a graph, often directed and/or time-varying.

%There is a vast literature in decentralized optimization, all of which we cannot cite in here.
%One of the well-studied class of algorithms is so called
The literature on distributed optimization methods is vast and involves first-order methods in the primal domain, the dual domain, augmented Lagrangian methods, or Newton methods, to name a few. Here we discuss methods that are most closely related to the method under consideration. Among those, one of the most well-studied techniques are the so called \textit{consensus-based} optimization algorithms \cite{Nedic2009,Nedic2010,Bianchi2012,tsianos2012consensus,RamNV09,RamNV10,Nedic11,KS2011,LobelOF2011,Lobel2011,ANAO,tsianos-pushsum,Duchi-dual,Wei-admm,Bianchi-Explicit,Bianchi-coord,Bianchi-asyn,WYin-Extra
} (see also the literature for the consensus problem itself \cite{con01,con03,con05}),
where the goal is to repeatedly average the estimates of all agents in a decentralized fashion
in order to obtain a network-wide consensus.
Between the averaging steps, each agent usually performs a single local optimization step.
Overall, the agents use their local information to cooperatively steer the consensus point toward the optimal set of the global problem.

Decentralized algorithms that fall in this class of methods can be distinguished based on which averaging scheme or optimization method is used, and in which space (primal or dual) the iterates are maintained.
All these algorithms often require
expensive optimization steps or exact projections on a complicated constraint set at every iteration.
%lots of computations like optimization in inner loops or projections on complicated constraint set can be made at every iteration.
Such intensive computations, however, require time and may shorten the lifespan of certain systems, such as wireless sensor networks or robotic networks.

In this work,
%we are specifically interested in distributed problems involving a huge number of complex inequalities,
%while the number of decision variables is moderate.
we propose a new approximate projection based decentralized algorithm %in primal space
and prove its convergence.
%Note that dual methods cannot be applied in this case due to large number of dual variables.
Our work in this paper is an extension of the author's previous work \cite{SL12b}.
Specifically, we use the same local information exchange model and gradient descent algorithm as in \cite{SL12b}, but a different projection method motivated by the work in \cite{Polyak2001409}.
In contrast to \cite{SL12b}, our contribution can be summarized as follows:
(1) Instead of using the Euclidean projection, we approximate it by measuring the constraint violation
and taking a subgradient step minimizing this violation;
(2) We show convergence under milder assumptions. Specifically, we remove the smoothness assumption in the objective functions;
(3) We propose two different averaging schemes to mitigate the disagreements among the agents' local estimates, one of which can lift the doubly stochasticity assumption on the weight matrices.
%We show that this decentralized algorithm terminates within a finite number of iterations when solving feasibility problems.

Considering that projections have a closed form solution only in a few special cases of constraints,
%only a few special cases of constraints have a closed form solution to the projection,
our new algorithm is more general and can be applied to a wider class of problems including Semidefinite Programming (SDP), where the constraints are represented by
Linear Matrix Inequalities (LMIs).
It is well known that even finding a feasible point that satisfies a handful of LMIs is a difficult problem on its own. The work in this paper is also related to the centralized
random projection algorithms for convex constrained optimization \cite{AN2011}
and convex feasibility problems \cite{nedicFP2011}.
Other related works are \cite{tempo2009,Calafiore:2010,Campi05uncertainconvex},
where optimization problems with uncertain constraints have been considered
by finding probabilistic feasible solutions through random sampling of constraints.

%Also, the computational overhead of this algorithm in each iteration is extremely cheap as the additional subgradient step requires only $O(n)$ calculations ($n$ is the number of decision variables).
% while, e.g.,an eigenvalue decomposition to exactly project onto semidefinite constraints requires $O(n^3)$ calculations.
%Therefore, our algorithm can also handle SDPs in a decentralized fashion with an extremely cheap per-iteration computational overhead.

The paper is organized as follows.
In Section \ref{sec:prob},
we formulate the optimization problem under consideration and discuss
specific problems of interest.
%we describe how the agents communicate, formalize the problem, and describe the specific class of problems of interest.
In Section \ref{sec:algo}, we provide our decentralized algorithm based on random approximate projections,
discuss the communication scheme employed by the agents,
state assumptions and the main results of this paper.
In Section \ref{sec:conv}, we first review some necessary results and lemmas from existing literature,
provide proofs of required lemmas, and then present the proofs of the main results discussed in Section \ref{sec:algo}.
In Section \ref{sec:sim}, we present simulation results for a decentralized SDP problem, which is optimal decentralized gossip averaging.
We conclude the paper with some comments in Section \ref{sec:con}.

\textbf{Notation: }
All vectors are viewed as column vectors. We write $x^{\top}$ to denote the transpose of a vector $x$.
The scalar product of two vectors $x$ and $y$ is $\la x, y\ra$.
For vectors associated with agent $i$ at time $k$, we use subscripts $i,k$ such as, for example, $p_{i,k}$, $x_{i,k}$, etc.
Unless otherwise stated, $\|\cdot\|$ represents the standard Euclidean norm.
For a set $S$, we use $|S|$ to denote its cardinality.
For a matrix $A \in \mathbb{R}^{n\times n}$, we use $[A]_{ij}$ to denote the entry of
the $i$-th row and $j$-th column and $\|A\|_F$ to denote the Frobenius norm $\|A\|_F = \left(\sum_{i,j=1}^n ([A]_{ij})^2\right)^{1/2}$. We use $\mathsf{Tr} A$ to denote the trace of $A$, i.e., $\mathsf{Tr} A = \sum_{i=1}^n [A]_{ii}$.
We denote by $\mathbb{S}_m$ and $\mathbb{S}_m^+$ the space of $m \times m$ real symmetric and real symmetric positive semidefinite matrices, respectively.
The matrix inequality $A \preceq 0$ means $-A$ is positive semidefinite.
We use $\1b$ and $\0b$ to denote vectors of all ones and zeros. The identity matrix is denoted by $I$.
We use $\textsf{Pr}\{Z\}$ and $\Es[Z]$ to denote the probability and the expectation of a random variable $Z$.
We write $\dist(x,\Xc)$ for the distance of a vector $x$ from a closed convex set $\Xc$, i.e., $\dist(x,\Xc) =
\min_{v\in \Xc} \|v-x\|$. We use $\proj_{\Xc}[x]$ for the Euclidean projection of a vector $x$ on the set $\Xc$, i.e., $\proj_{\Xc} [x] = \argmin_{v\in\Xc} \|v-x\|^2$.
We often abbreviate \textit{almost surely} and \textit{independent identically distributed} as \textit{a.s.} and \textit{i.i.d.}, respectively.

\section{Problem Definition\label{sec:prob}}
Consider a multiagent network system whose communication at time $k$ is governed by a digraph $\Gc_k = (\Vc,\Ec_k)$, where $\Vc = \{1,\ldots,N\}$ and $\Ec_k \subseteq \Vc \times \Vc$.
If there exists a directed link from agent $j$ to $i$, which we denote by $(j,i)$, agent $j$ may send its information to agent $i$. Thus, each agent $i \in \Vc$ can directly receive information only from the agents in its in-neighborhood
\begin{equation}\label{eqn:nbrin}
\Nc_{i,k}^{\text{in}} = \{j \in \Vc \mid (j,i) \in \Ec_k\} \cup \{i\},
\end{equation}
and send information only to the agents in its out-neighborhood
\begin{equation}\label{eqn:nbrout}
\Nc_{i,k}^{\text{out}} = \{j \in \Vc \mid (i,j) \in \Ec_k\} \cup \{i\},
\end{equation}
where in both $\Nc_{i,k}^{\text{in}}$ and $\Nc_{i,k}^{\text{out}}$,
we assume there exists a self-loop $(i,i)$ for all $i \in \Vc$.
Also, we use $\textrm{deg}_{i,k}$ to denote the in-degree of node $i$ at time $k$, i.e.,
\begin{equation}\label{eqn:indeg}
\textrm{deg}_{i,k} = |\Nc_{i,k}^{\text{in}}|.
\end{equation}

%$[W(k)]_{ij} = w_{ij}(k)$ is the weight assigned on the link $(j,i) \in \Ec(k)$.

\subsection{Problem Statement}
Our goal is to let the network of agents cooperatively solve the following convex minimization problem:
\begin{align}\label{eqn:prob}
\min_{x} ~& f(x) \triangleq \sum_{i\in \Vc} f_i(x) \\
\text{s.t. } & x \in \Xc,\quad \Xc \triangleq \Xc_0 \medcap \left(\medcap_{i\in \Vc} \Xc_i\right), \nonumber
%\text{s.t. } & x \in  \Xc_0 \medcap \left(\medcap_{i\in \Vc} \Xc_i\right), \nonumber
\end{align}
where only agent $i$ knows the function $f_i:\mathbb{R}^n \to \mathbb{R}$ and the constraint set $\Xc_i \subseteq \mathbb{R}^n$.
The set $\Xc_0 \subseteq \mathbb{R}^n$ is common to all agents and assumed to have some simple structure in the sense that the projection onto $\Xc_0$ can be made easily (e.g., a box, ball, probability simplex, or even $\mathbb{R}^n$).
% Note that the local constraints may overlap, i.e., $\Xc_i \cap \Xc_j$ for $i \neq j$ may be nonempty. Also,
Note that the common constraint set, i.e., $\Xc_0 = \Xc_i$ for all $i \in \Vc$, is a special case of this problem definition. We assume that the set of optimal solutions $\Xc^* = \argmin_{x \in \Xc} f(x)$ is nonempty.

\sml{
We assume each agent $i$'s local constraint set $\Xc_i$ consists of one or more algebraic inequalities, which we denote by
\[
\Xc_i = \{x \in \Rbb^n \mid g(x,\o) \le 0,~\forall \o \in \O_i\},
\]
where $\O_i$ is a finite collection of indices.
%and $|\O_i|$ is the number of inequalities (or constraint components) in agent $i$'s local set $\Xc_i$.
From this definition, the feasible set $\Xc$ can be precisely represented as
\[
\Xc = \{x \in \Xc_0 \mid g(x,\o) \le 0,~\forall \o \in \O_i,i\in \Vc\}.
\]
Note that some of the inequalities may overlap across different agents, i.e., $\O_i \cap \O_j$ for $i \neq j$ can be either empty or nonempty.}

We also consider an equivalent epigraph form of problem \eqref{eqn:prob} by introducing a new set of variables $t = [t_1 ~\ldots~ t_N]^{\top}\in \mathbb{R}^N$. Consider that the local constraint set $\Xc_i$ for $i\in \Vc$ now includes the additional inequality constraint
$f_i(y) \le t_i$. Then, problem \eqref{eqn:prob} is equivalent to:
\begin{align}\label{eqn:probepi}
\min_{x} ~&  a^{\top} x \\
\text{s.t. } & y \in \Xc,\quad \Xc \triangleq \Xc_0 \medcap \left(\medcap_{i\in \Vc} \Xc_i\right), \nonumber
%\text{s.t. } & x \in  \Xc_0 \medcap \left(\medcap_{i\in \Vc} \Xc_i\right), \nonumber
\end{align}
where $x = [y^{\top}~ t^{\top}]^{\top}$ and $a = [\0b^{\top}~\1b^{\top}]^{\top}$.
Note that this reformulation has also been introduced in \cite{Keyou}
for solving distributed robust convex optimization.
In the case of a single agent, the algorithm can be also seen as scenario approach for solving robust convex optimization (see e.g., \cite{SA1,SA2,SA3} and references therein).

\subsection{Problems of Interest\label{ssec:poi}}
\sml{Problems of particular interest are those involving lots of nontrivial constraints on which exact projections are impossible or computationally intractable.
%Therefore, when nodes are limited with computational capability, we can certainly benefit from a decentralized framework.
Here we provide two such examples:}
\begin{enumerate}
\item Robust Linear Inequalities:
\begin{align*}
\Xc_i & =  \Big\{x \in \mathbb{R}^n \mid  A(\o)x \le b(\o),~ \sml{\forall \o \text{ such that}}\\
&\sml{\|A(\o)-A_0\|_{\text{op}} \le r_1
 \text{ and } \|b(\o)-b_0\|_{\text{op}} \le r_2}\Big\},
\end{align*}
where $A_0 \in \mathbb{R}^{m\times n}, b_0 \in\mathbb{R}^m$ are nominal data, $r_1, r_2 \ge 0$ are the levels of uncertainty, and $\|\cdot\|_{\text{op}}$ denotes an operator norm. Here we can not handle each row of $A(\o)x \le b(\o)$ separately as in \cite{SL12b} due to the matrix operator norm $\|\cdot\|_{\text{op}}$ .
%If each row of $A(\o_i)x \le b(\o_i)$ can be handled separately, we can use, for example, random projection techniques in \cite{SL12b}.
\item Linear Matrix Inequalities:
\begin{align}\label{eqn:LMI}
%\Xc_i = \Big\{x \in \mathbb{R}^n \mid & A_0(\o) + \sum_{j=1}^n x_j A_j(\o) \preceq 0, \nonumber\\
%& \text{for } \o \in \O_i \Big\},
\Xc_i = \Big\{x \in \mathbb{R}^n \mid & A_0(\o) + \sum_{j=1}^n x_j A_j(\o) \preceq 0,~\sml{\forall \o \in \O_i} \Big\},
\end{align}
where $A_j(\o) \in \mathbb{S}_m$ for $j = 0, 1,\ldots,n$, $\o\in \O_i$ are given matrices.
The inequalities in (\ref{eqn:LMI}) are referred to as \textit{linear matrix inequalities} (LMIs).
A semidefinite programming (SDP) problem has one or more LMI constraints.
Finding a feasible point of the set (\ref{eqn:LMI}) is often a difficult problem on its own.
\end{enumerate}
Note that the inequalities in (\ref{eqn:LMI}) can represent a wide variety of convex constraints (see \cite{BEFB:94} for more details).
For example, quadratic inequalities, inequalities involving matrix norms, and various inequality constraints arising in robust control such as Lyapunov and quadratic matrix inequalities can be all cast as LMIs in (\ref{eqn:LMI}).
When all matrices $A_j(\o)$ in (\ref{eqn:LMI}) are diagonal, \sml{the LMIs reduce to regular linear inequalities}.

%Let $M_i$'s are disjoint (i.e., $M_i \cup M_j = \emptyset$ if $i \neq j$ and $M = \cup_{i=1}^N M_i$.
%Particular cases of problem are:

\section{Algorithm, Assumptions, and Main Results\label{sec:algo}}
Our goal is to design a decentralized protocol by which each agent $i \in \Vc$ maintains a sequence of the local copy $\{x_{i,k}\}_{k \ge 0}$ converging to the same point in $\Xc^*$ as $k$ goes to infinity.
Since we assume that the local constraint sets $\Xc_i$'s are nontrivial, we do not find an exact projection onto $\Xc_i$ at each step of the algorithm. Instead, at iteration $k$, each agent $i$ randomly generates an index $\o_{i,k} \in \O_i$ and makes an \textit{approximate projection} on the selected inequality $g(\cdot,\o_{i,k}) \le 0$.

\subsection{Decentralized Algorithm with Approximate Projections}
We formally present our decentralized algorithm, named the Decentralized Approximate Projection (DAP), in Algorithm \ref{alg:DAP}.
Each agent $i$ maintains a sequence $\{x_{i,k}\}_{k \ge 0}$. The element $x_{i,k}$ of the sequence can be seen as the agent $i$'s estimate of the decision variable $x$ at time $k$.
Let $g^+(x,\o)$ denote the function that measures the violation of the constraint $g(\cdot,\o)$ at $x$, i.e.,
$g^+(x,\o) \triangleq \max\{g(x,\o),0\}$.

\begin{algorithm}[t]
\caption{Decentralized Approximate Projection (DAP)}
\label{alg:DAP}
\begin{algorithmic}
\STATE{Let $x_{i,0} \in \Xc_0$ for $i \in \Vc$ and the nonnegative parameter $\{\a_k\}_{k\ge0}$ be given.}
\STATE{Set $k:=0$}
\WHILE{Maximum iteration number is reached}
\STATE{Each agent $i$ updates $x_{i,k}$ according to}
\begin{subequations}
\begin{align}
p_{i,k} =~&  \sum_{j\in \Vc} [W_k]_{ij} x_{j,k-1}\label{eqn:algo1}\\
v_{i,k} =~& \mathsf{\Pi}_{\Xc_0}\left[p_{i,k} - \a_k s_{i,k}\right]\label{eqn:algo2}\\
%x_{i,k} =~& \mathsf{\Pi}_{\Xc_0}\left[v_{i,k} - \frac{g^+_i(v_{i,k},\o_{i,k})}{\|d_{i,k}\|^2}d_{i,k} \right], \label{eqn:algo3}
x_{i,k} =~& \mathsf{\Pi}_{\Xc_0}\left[v_{i,k} - \frac{g^+(v_{i,k},\o_{i,k})}{\|d_{i,k}\|^2}d_{i,k} \right]
\label{eqn:algo3}
\end{align}
\end{subequations}
\STATE{Set $k:=k+1$}
\ENDWHILE
%\RETURN $\bar{\zb}^N = (\sum_{k=1}^N \g_k)^{-1} \sum_{k=1}^N \g_k \zb^{k}$.
\end{algorithmic}
\end{algorithm}

At $k = 0$, the estimates $x_{i,0}$ are locally initialized such that $x_{i,0} \in \Xc_0$.
At time step $k$, all agents $j \in \Vc$ broadcast their previous estimates $x_{j,k-1}$ to all of the nodes in their out-neighborhood, i.e., to all agents $i$ such that $(i,j) \in \Ec_k$.   Then, each agent $i \in \Vc$ updates $x_{i,k}$ using \eqref{eqn:algo1}-\eqref{eqn:algo3},
where $W_k$ is a nonnegative $N\times N$ weight matrix, $\{\a_k\}$ is a positive sequence of nonincreasing stepsizes;
$s_{i,k}$ is a subgradient of the function $f_i$ at $p_{i,k}$; $\o_{i,k}$ is a random variable taking values in the index set $\O_i$; and $d_{i,k}$ is a subgradient of $g^+(\cdot,\o_{i,k})$ evaluated at $v_{i,k}$.
The vector $d_{i,k}$ is chosen such that $d_{i,k} \in \partial g^+(v_{i,k},\o_{i,k})$ if $g^+(v_{i,k},\o_{i,k}) >0$, and $d_{i,k} = d_i$ for some $d_i \neq 0$, if $g^+(v_{i,k},\o_{i,k}) =0$.
%The scalar $\lambda_{i,k}$ is a nonnegative stepsize which will be defined later.
%that is proportional to the constraint violation $g^+(v_{i,k},\o_{i,k})$ if $g^+(v_{i,k},\o_{i,k}) >0$;

More specifically, in (\ref{eqn:algo1}), each agent $i$ calculates a weighted average of the received messages (including its own message $x_{i,k-1}$) to obtain $p_{i,k}$. Specifically, $[W_k]_{ij} \ge 0$ is the weight that agent $i$ allocates to the message $x_{j,k-1}$. This communication step is decentralized since the weight matrix $W_k$ respects the topology of the graph $\Gc_k$, i.e., $[W_k]_{ij} >0$ only if $(j,i) \in \Ec_k$ and $[W_k]_{ij} =0$, otherwise.
In (\ref{eqn:algo2}), each agent $i$ adjusts the average $p_{i,k}$ in the direction of the negative subgradient  ($-s_{i,k}$)  of its local objective $f_i$ to obtain $v_{i,k}$. The adjusted average is projected back to the simple set $\Xc_0$.
In (\ref{eqn:algo3}), agent $i$ observes a random realization of $\o_{i,k} \in \O_i$ and measures the feasibility violation of the selected component constraint $g(\cdot,\o_{i,k})$ at $v_{i,k}$. If $g^+(v_{i,k},\o_{i,k}) > 0$, it calculates a subgradient $d_{i,k} \in \partial g^+(v_{i,k},\o_{i,k})$ and takes an additional subgradient step with the stepsize $\frac{g^+(v_{i,k},\o_{i,k})}{\|d_{i,k}\|^2}$ to minimize this violation.
\sml{If $g^+(v_{i,k},\omega_{i,k}) = 0$, then the current point $v_{i,k}$ already satisfies the selected inequality $g(\cdot,\omega_{i,k})\le 0$. In this case, there is no need to move the point further into the selected set. Therefore, the approximate projection step \eqref{eqn:algo3} is just omitted.}
%The stepsize rule is proportional to the amount of violation $g_i^+(v_{i,k},\o_{i,k})$ as dividing by $\|d_{i,k}\|^2$ works like a normalization of the direction $d_{i,k}$.
%The adjusted vector is then projected back to the set $\Xc_0$.

Note that the description of the DAP algorithm is only conceptual at this moment
since we have not specified the parameters $\{\a_k\}$ and $\{W_k\}$ yet.
The stepsizes $\{\a_k\}$  should be nonnegative, nonincreasing and such that
\begin{align}\label{eqn:step}
\sum_{k=1}^{\infty}\a_k = \infty \text{ and } \sum_{k=1}^{\infty}\a_k^2 < \infty.
\end{align}
For the sequence of weight matrices $\{W_k\}$, we assume the following.

\begin{assumption}\label{assume:ds}
\textit{For all $k \ge 1$,
\begin{enumerate}
\item[(a)] $[W_k]_{ij} \ge 0$ for all $i,j \in \Vc$ and $[W_k]_{ij} = 0$ only if $j \not\in \Nc_{i,k}^{\text{in}}$.
\item[(b)] There exists a scalar $\nu \in (0,1)$ such that $[W_k]_{ij} \ge \nu$ only if $j \in \Nc_{i,k}^{\text{in}}$.
\item[(c)] $\sum_{j\in \Vc} [W_k]_{ij} = 1$ for all $i \in \Vc$, and $\sum_{i\in \Vc} [W_k]_{ij} = 1$ for all $j \in \Vc$.
\end{enumerate}
}
\end{assumption}
Condition (a) ensures that the weight matrices $W_k$ respect the underlying topology $\Gc_k$ for every $k$ so that the communication is indeed decentralized.
The lower boundedness of the weights in (b)
is required to show consensus among all agents (see \cite{ANquan2008} for more details)
but the agents need not know the $\nu$ value in running the algorithm.
Condition (a) and (c) imply doubly stochasticity of the matrices $W_k$. %is required for balanced communication.

Note that assuming doubly stochasticity of $W_k$ for all $k \ge 1$ often requires a global view of the network (unless the underlying graph is regular or fully connected) and
not all directed graphs admit a doubly stochastic matrix \cite{GharesifardCortes,AsuCDC2015}.
The linear objective function in the epigraph formulation \eqref{eqn:probepi} allows us to lift this assumption
and to use a weight-imbalanced mixing matrix.
For example, we can assume a time-invariant network, i.e., $\Gc_k = \Gc$ for $k \ge 1$, and employ a row stochastic matrix $W$ whose entries are defined as follows:
\begin{align}\label{eqn:wgtrow}
[W]_{ij} =
\left\{
\begin{array}{ll}
\frac{1}{\textrm{deg}_i} & \text{ if } j \in \Nc_i^{\text{in}},\\
0 & \text{ otherwise.}
\end{array}
\right.
\end{align}
%where $\textrm{deg}_i$ is defined in \eqref{eqn:indeg}.
Note that this choice of weights also respects the underlying topology $\Gc$. Moreover,
the matrix $W$ is not necessarily column stochastic.

\sml{
Recall that our problems of interest involve a large number of constraints.
Therefore, the random selection of a constraint in \eqref{eqn:algo3}
serves as a computationally efficient alternative
to finding the most violated constraint,
which typically has significantly higher per-iteration complexity.
%as  may be a hard task.
%However, note that we show in Proposition \ref{prop:feasibility} that
%Selecting the most violated constraint at each iteration can be beneficial, i.e., in guaranteeing finite termination of the algorithm for distributed feasibility problems, as we show in Proposition \ref{prop:as2}. (TODO)
%at the cost of relatively high computational complexity.
%as we assume
%An alternative would be
%to measure constraint violations for all constraint components and selecting the most violated one
%would require lots of computations per iteration.
Another situation that necessitates the random selection approach is when the constraints are not fully given in advance, but are rather revealed in a sequential fashion (as in online optimization).}

Note that the step \eqref{eqn:algo3} guarantees that $x_{i,k} \in \Xc_0$ for all $k \ge 0$ and $i \in \Vc$, but it does not necessarily guarantee $x_{i,k}\in \Xc$. Nevertheless, in Section \ref{sec:conv}, we show that $x_{i,k}$ for all $i \in \Vc$ asymptotically achieve feasibility, i.e., $\lim_{k\to \infty}\|x_{i,k}-\proj_{\Xc}[p_{i,k}]\| = 0$ for all $i \in \Vc$.

%\begin{comment}
To further explain the step (\ref{eqn:algo3}), let us consider the two particular cases mentioned in Section \ref{ssec:poi}.
%Note that exact projections onto the sets described by the following two cases are intractable.
\begin{enumerate}
\item Let $c^+$ denote a projection of a vector $c \in \mathbb{R}^m$ onto the nonnegative orthant.
We introduce a scalar function in order to handle all the rows of the inequality $A(\o)x \le b$ concurrently,
\[
g^+(x,\o) = \|(A(\o)x-b)^+\|
\]
which is convex in $x$ for any given $\o \in \O_i$ and $i \in \Vc$ \cite[Chapter 3.2]{Boyd}. Then, it is straightforward to see
that its subgradient can be calculated as
\[
\partial g^+(x,\o) = \frac{A(\o)^{\top}(A(\o)x-b)^+}{\|(A(\o)x-b)^+\|},
\]
if $g^+(x,\o)>0$, and we can use $\partial g^+(x,\o) = 0$, otherwise.
\item Let us define the projection $A^+$ of a real symmetric matrix $A$ onto the cone of positive semidefinite matrices.
For any $A \in \mathbb{S}^m$, we can find an eigenvalue decomposition
$
A = B\Lambda B^{\top},
$
where $B$ is an orthogonal matrix and $\Lambda = \diag(\lambda_1,\ldots,\lambda_m)$. Then, its projection is given by
\[
A^+ = B\Lambda^+B^{\top},
\]
where $\Lambda^+ = \diag(\lambda_1^+,\ldots,\lambda_m^+)$ with $\lambda_i^+ = \max\{0,\lambda_i\}$ \cite{Polyak-grad}.
Let us define
\[
A(x,\o) \triangleq A_0(\o) + \sum_{j=1}^n x_j A_j(\o).
\]
Then, the amount of violation of the corresponding LMI constraint $A(x,\o) \preceq 0$ can be measured by the following convex scalar function:
\begin{align}\label{eqn:Fro}
g^+(x,\o) = \|A^+(x,\o)\|_F.
\end{align}
By direct calculations, it is not difficult to see that its subgradient is given by
\begin{equation}\label{eqn:subgradLMI}
\partial g^+(x,\o) = \frac{1}{g^+(x,\o)}
\left(\begin{array}{c}
\mathsf{Tr} A_1 A^+(x,\o)\\
\vdots\\
\mathsf{Tr} A_n A^+(x,\o)
\end{array}\right),
\end{equation}
if $g^+(x,\o) > 0$, and we can use $\partial g^+(x,\o) = 0$, otherwise \cite[Lemma 1]{Polyak2001409}.
\end{enumerate}
%\end{comment}

\sml{
\begin{remark}
Note that the computational complexity of step \eqref{eqn:algo3} depends on the type of the function $g(\cdot,\omega_{i,k})$. If $g(\cdot,\omega_{i,k})$ is a general convex function, it takes $O(1)$ computations for the evaluation of $g^+(\cdot,\omega_{i,k})$ and $O(n)$ computations for the evaluation of the gradient $d_{i,k}$.
If $g(\cdot,\omega_{i,k})$ is an LMI constraint, it takes $O(m^3)$ computations in the worst-case for the eigenvalue decomposition and $O(m^2)$ for the computation of the Frobenius norm (cf. Eq. \eqref{eqn:Fro}). We would also need $O(m^2n)$ computations for computing the traces (cf. Eq. \eqref{eqn:subgradLMI}).
This eigenvalue decomposition is necessary for
projection (or approximate projection) onto the cone of positive semidefinite matrices. $\blacksquare$
\end{remark}}
%prevailing in SDP algorithms. Indeed, it is necessary for measuring the feasibility violation (for identifying the positive eigenvalues) and for projection onto the cone of positive semidefinite matrices (for setting the negative eigenvalues zero).

It is also worth mentioning that the algorithm (\ref{eqn:algo1})-(\ref{eqn:algo3})
includes the method that has been proposed in \cite{SL12b} as a special case. In order to see this,
let $\Xc_0 = \mathbb{R}^n$ and $g(x,\o_{i,k}) = \dist(x,\Xc_i^{\o_{i,k}})$,
where $\Xc_i^{\o_{i,k}} = \{x \in \mathbb{R}^n \mid g(x,\o_{i,k}) \le 0\}$.
%Suppose we use $\lambda_{i,k}$ in \eqref{eqn:lda} with $r=0$.
Then, it is not difficult to see that
\[
d_{i,k} = \frac{v_{i,k}-\mathsf{\Pi}_{\Xc_i^{\o_{i,k}}}[v_{i,k}]}{\|v_{i,k}-\mathsf{\Pi}_{\Xc_i^{\o_{i,k}}}[v_{i,k}]\|},
\]
and since $\dist(v_{i,k},\mathsf{\Pi}_{\Xc_i^{\o_{i,k}}}) = \|v_{i,k}-\mathsf{\Pi}_{\Xc_i^{\o_{i,k}}}[v_{i,k}]\|$,
the steps (\ref{eqn:algo2})-(\ref{eqn:algo3}) reduce to
\[
x_{i,k} = \mathsf{\Pi}_{\Xc_i^{\o_{i,k}}}[p_{i,k}-\a_ks_{i,k}],
\]
which is exactly the algorithm in \cite{SL12b}.

\subsection{Assumptions}
For the optimization problem (\ref{eqn:prob}), we make the following assumptions on the set $\Xc_0$, the objective functions $f_i(x)$ for $i\in \Vc$, and the constraint functions $g(x,\o)$ for $\o \in \O_i$ and $i\in \Vc$.

\begin{assumption}\label{assume:f}
\textit{
\begin{enumerate}
\item[(a)] The set $\Xc_0$ is nonempty, closed and convex.
\item[(b)] The function $f_i(x)$, for each $i \in \Vc$, is defined and convex (not necessarily differentiable) over some open set that contains $\Xc_0$.
\item[(c)] The subgradients $s \in \partial f_i (x)$ are uniformly bounded over the set $\Xc_0$. That is, for all $i \in \Vc$, there is a scalar $C_{f_i}$ such that for all $s \in \partial f_i (x)$ and  $x \in \Xc_0$,
    \[
    \|s\| \le C_{f_i}.
    \]
\item[(d)] \sml{The function $g(x,\o)$, for each $\o \in \O_i$ and $i\in \Vc$, is defined and convex in $x$ (not necessarily differentiable) over some open set that contains $\Xc_0$.}
\item[(e)] The subgradients $d \in \partial g^+ (x,\o)$ are uniformly bounded over the set $\Xc_0$. That is, there is a scalar $C_{g}$ such that for all $d\in \partial g^+ (x,\o),~ x \in \Xc_0,~ \o \in \O_i$, and $i \in \Vc$,
    \[
    \|d\| \le C_{g}.
    \]
\end{enumerate}
}
\end{assumption}
By Assumption \ref{assume:f}, the subdifferentials $\partial f_i (x)$ and $\partial g^+ (x,\o)$ are nonempty over $\Xc_0$. It also implies that for any $i \in \Vc$ and $x, y \in \Xc_0$,
\begin{equation}\label{eqn:Cfi}
|f_i(x)-f_i(y)| \le C_{f_i} \|x-y\|,
\end{equation}
and for any $\o \in \O_i$, $i \in \Vc$, and $x, y \in \Xc_0$,
\begin{equation}\label{eqn:Cgi}
|g^+(x,\o)-g^+(y,\o)| \le C_{g} \|x-y\|.
\end{equation}
One sufficient condition for Assumption \ref{assume:f}(c) and \ref{assume:f}(e) is that the set $\Xc_0$ is compact.

%Let $\Fc_k$ be the filteration of the algorithm's history until time $k$.
We also require the following two assumptions.
\begin{assumption}\label{assume:Pr}
\textit{We assume that $\o_{i,k} \in \O_i$ are \textit{i.i.d.} samples from some probability distribution on $\O_i$
and independent across agents.
Furthermore, each $\O_i$ is a finite set and each element of $\O_i$ is generated with nonzero probability,
%Any random variable $\o \in \O_i$ has a positive probability distribution over the set $\O_i$,
i.e., for any $\o \in \O_i$ and $i \in \Vc$}
\[
\textsf{Pr}\{\o \mid \o\in\O_i\} > 0
\]
%If for some $x \in \Xc_0$, $x \not\in \Xc_i$, then the probability of generating a violated inequality is nonzero, i.e.,
%for any $i \in \Vc$
%\[
%\Pr\{\o \mid g^+(x,\o) > 0 \text{ and } \o \in \O_i\} > 0
%\]
\end{assumption}

Let $\Fc_k$ denote the algorithm's history up to time $k$. i.e.,
\[
\Fc_k = \{x_{i,0}, ~(\o_{i,t},~1\le t \le k),~ i \in \Vc\},
\]
and $\Fc_= \{x_{i,0},~i\in \Vc\}$.

\begin{assumption}\label{assume:c}
\textit{For all $i \in \Vc$,
%For all $i \in \Vc$ and for any random variable $\o \in \O_i$ with a positive probability distribution over the set $\O_i$,
there exists a constant $c > 0$ such that for all $x \in \Xc_0$
\[
\dist^2 (x,\Xc) \le c \Es\left[(g^+(x,\o))^2 \mid \Fc_{k-1}\right],
\]
where the expectation is taken with respect to the set $\Omega_i$.
}
\end{assumption}
The upper bound in Assumption \ref{assume:c} is known as \textit{global error bound} and is crucial for the convergence analysis of our method (\ref{eqn:algo1})-(\ref{eqn:algo3}).
Sufficient conditions for this bound have been shown in \cite{pang-book} and \cite{Lewis96errorbounds}, which require the existence of a Slater point, i.e.,
let $\Xc_0 = \{x\mid g_0(x) \le 0\}$, then
there exists a point $\tilde{x}$ such that $g_0(\tilde{x}) < 0$ and $g(\tilde{x},\o) < 0$ for all $\o$.
When each function $g(\cdot,\o)$ and $g_0(\cdot)$ is either a linear equality or inequality,
Assumption \ref{assume:c} is called \textit{linear regularity} and can be shown to hold by using the results in \cite{Bauschke:1996} and \cite{Burke93weaksharp} (see also \cite{Deutsch200636,Deutsch200656,Deutsch2008155}).
\begin{comment}
Also, if the random sequences $\{\o_{i,k}\}$ are stationary Markov processes (which are independent of the distribution of the initial points $x_{i,0}$),
then Assumption \ref{assume:c} holds when the feasible set $\Xc$ has a nonempty interior \cite{Gubin19671}.
\end{comment}
\begin{comment}
The lower bound in Assumption \ref{assume:c} is is usually not difficult to show. Indeed, when $g(x,\o_{i,k}) = \dist(x, \Xc_i^{\o_{i,k}})$, where $\Xc_i^{\o_{i,k}} = \{x \mid g(x,\o_{i,k}) \le 0\}$ (in this case the step (\ref{eqn:algo3}) just becomes an exact projection) the bound holds immediately with $c_1 = 1$. However, we write this bound as an assumption as the constant $c_1$ will be varying for different kind of constraint functions and we want to prove the convergence of the algorithm for general constraint functions.
\end{comment}
%Assumption \ref{assume:Pr} is satisfied, for example, when each $\O_i$ is a finite set and each element of $\O_i$ is generated with nonzero probability. (TODO: explanation)

The inter-agent communication relies on the time-varying graph sequence $\Gc_k = (\Vc,\Ec_k)$, for $k\ge 1$. A key assumptions on these communication graphs is the following:
\begin{assumption}\label{assume:nc}
\textit{
There exists a scalar $Q$ such that the graphs
$\left(\Vc,\bigcup_{\ell = 0,\ldots,Q-1} \Ec_{k+\ell}\right)$
are strongly connected for all $k\ge 1$.
}
\end{assumption}

\sml{Assumption~\ref{assume:nc} ensures that there exists a path
from one agent to every other agent within any bounded interval of length $Q$.
We say that such a sequence of graphs is $Q$-strongly connected.}

%sufficiently often so that all functions and all constraints ($f_i$'s and $\mathcal{X}_i$'s)
%influence the iterates of all agents.

\subsection{Main Results}

%Our first proposition demonstrates the correctness of the algorithm (\ref{eqn:algo1})-(\ref{eqn:algo3}). The end result is stated in the following proposition, which holds under the assumptions we have laid out above.

We now provide two convergence results for DAP, Proposition \ref{prop:as} for problem \eqref{eqn:prob} and Proposition \ref{prop:as2} for problem \eqref{eqn:probepi}, respectively, for which we use different assumptions on the network.
The first proposition states a convergence result which holds under a $Q$-strongly connected time-varying sequence of graphs $\{\Gc_k\}$  and the corresponding doubly stochastic matrices $\{W_k\}$ which respect the graph topologies.
\begin{proposition}\label{prop:as}
\textit{
Let Assumptions \ref{assume:ds} - \ref{assume:nc} hold and the optimal set $\Xc^*$ of \eqref{eqn:prob} be nonempty.
Let the nonnegative and nonincreasing stepsize $\{\a_k\}$ satisfy conditions in \eqref{eqn:step}.
Then, %with the choice of $r = 0$ in \eqref{eqn:lda},
the iterates $\{x_{i,k}\}$ generated by each agent $i \in \Vc$ via DAP in Algorithm \ref{alg:DAP}
%with the choice of weight in \textbf{Option 1}
converge almost surely to the same point in the optimal set $\Xc^*$ of \eqref{eqn:prob}, i.e., for a random point $x^{\star}\in \Xc^*$
\[
\lim_{k\to\infty} x_{i,k} = x^{\star} \quad \text{for all } i \in \Vc \quad a.s.
\]
}
\end{proposition}

The second proposition states a convergence result which holds under a connected time-invariant graph $\Gc$ and the corresponding row stochastic matrix $W$ which respects the graph topology. Note that for this case we can set $f_i(x) = \frac{1}{|\Vc|}a^{\top}x$ for all $i \in \Vc$ and $s_{i,k} = a$ for all $i \in \Vc$ and $k \ge 0$.
\begin{proposition}\label{prop:as2}
\textit{
Let Assumptions \ref{assume:f} - \ref{assume:c} hold and the optimal set $\Xc^*$ of \eqref{eqn:probepi} be nonempty.
Let the nonnegative and nonincreasing stepsize $\{\a_k\}$ satisfy conditions in \eqref{eqn:step}.
Let $W_k = W$ for all $k\ge 1$ and
Assumption \ref{assume:nc} hold with $Q=1$.
Then, %with the choice of $r = 0$ in \eqref{eqn:lda},
the iterates $\{x_{i,k}\}$ generated by each agent $i \in \Vc$ via DAP in Algorithm \ref{alg:DAP}
with the choice of weight in \eqref{eqn:wgtrow}
converge almost surely to the same point in the optimal set $\Xc^*$ of \eqref{eqn:probepi}, i.e., for a random point $x^{\star}\in \Xc^*$
\[
\lim_{k\to\infty} x_{i,k} = x^{\star} \quad \text{for all } i \in \Vc \quad a.s.
\]
}
\end{proposition}

In the next paragraph, we provide some intuition on why DAP still converges even with a choice of a row stochastic matrix $W$.
Note that the use of row stochastic matrix results in ``biased'' consensus, which is related to the left-eigenvector, see e.g., \cite{Meyer:2000}. The following lemma states this well-known result.
\begin{lemma}\label{lem:perron}
Let Assumption \ref{assume:nc} hold with $Q=1$ and $W = W_k$ for any $k \ge 1$. Then, there exists a normalized left-eigenvector $\pi \in \mathbb{R}^N$ such that
\begin{align*}
\pi^{\top}W = \pi^{\top}.
\end{align*}
Moreover, $[\pi]_i>0$ for all $i \in \Vc$.
\end{lemma}

This will lead
the algorithm to converge to an optimal solution of a biased objective function, $\sum_{i\in\Vc} [\pi]_i f_i(x)$,
instead of the true objective $\sum_{i\in\Vc} f_i(x)$.
However, since the $f_i(x)$'s in the reformulated epigraph form in \eqref{eqn:probepi} is now linear,
the biased objective function
is just the same as the original function, i.e.,
$\sum_{i \in \Vc}[\pi]_i f_i(x) = \sum_{i \in \Vc}[\pi]_i a^{\top}x = a^{\top}x$.

\section{Convergence Analysis\label{sec:conv}}
In this section, we are concerned with demonstrating the convergence results stated in Proposition \ref{prop:as} and \ref{prop:as2}. First we review some lemmas from existing literature that are necessary in our analysis.

\subsection{Preliminary Results}
First we state a non-expansiveness property of the projection operator (see \cite{BNO} for its proof).
\begin{lemma}\label{lem:proj}
\textit{
Let $\mathcal{X} \subseteq \mathbb{R}^n$ be a nonempty closed convex set.
The function $\mathsf{\Pi}_{\mathcal{X}}: \mathbb{R}^n \rightarrow \mathcal{X}$ is nonexpansive, i.e.,
\[
\|\mathsf{\Pi}_{\mathcal{X}}[x]-\mathsf{\Pi}_{\mathcal{X}}[y]\| \leq \|x-y\|
\quad \hbox{for all }x,y \in \mathbb{R}^n.
\]
}
\end{lemma}

In our analysis of the algorithm, we also make use of the following convergence result due to Robbins and Siegmund (see~\cite[Lemma 10-11, p. 49-50]{polyak}).
\begin{theorem}\label{thm:super}
\textit{
Let $\{v_k\}$, $\{u_k\}$, $\{a_k\}$ and $\{b_k\}$ be sequences of non-negative random variables such that
\[
\mathsf{E}[v_{k+1}|\Fc_k] \leq (1+a_k)v_k - u_k + b_k
\quad\text{ for all } k \geq 0 \quad a.s.
\]
where $\Fc_k$ denotes the collection $v_0,\ldots,v_k$, $u_0,\ldots,u_k$, $a_0,\ldots,a_k$ and $b_0,\ldots,b_k$.
Also, let $\sum_{k=0}^{\infty} a_k < \infty$ and $\sum_{k=0}^{\infty} b_k < \infty$ \textit{a.s.} Then, we have $\lim_{k \rightarrow \infty} v_k = v$ for a random variable $v \geq 0$ \textit{a.s.}, and $\sum_{k=0}^{\infty} u_k < \infty$ \textit{a.s.}
}
\end{theorem}

%\noindent\textit{Convexity and Doubly Stochasticity:~}
In the following lemma, we show a relation of $p_{i,k}$ and $x_{i,k-1}$
associated with any convex function $h$ %and any arbitrary vector $\check{x} \in \Xc_0$
which will be often used in the analysis. For example, $h(x) = \|x-a\|^2$ for some $a\in\mathbb{R}^n$ or
$h(x) = \dist^2 (x,\Xc)$.
\begin{lemma}\label{lem:ds}
\textit{
Let Assumption \ref{assume:ds} hold.
Then, for any convex function $h : \mathbb{R}^n \to \mathbb{R}$, we have
\[
\sum_{i\in \Vc} h(p_{i,k}) \le  \sum_{i\in \Vc} h(x_{i,k-1})
\]
}
\end{lemma}
\begin{IEEEproof}
The doubly stochasticity of the weights plays a crucial role in this lemma. From the definition of $p_{i,k}$ in (\ref{eqn:algo1}),
\begin{align*}
\sum_{i\in \Vc} h(p_{i,k}) \le~& \sum_{i\in \Vc} \sum_{j\in \Vc} [W_k]_{ij} h(x_{j,k-1})\\
=~& \sum_{j\in \Vc} \left(\sum_{i\in \Vc} [W_k]_{ij}\right)h(x_{j,k-1})\\
=~& \sum_{j\in \Vc} h(x_{j,k-1}).
\end{align*}
\end{IEEEproof}

%\noindent\textit{Consensus with errors:~}
Lastly, for the convergence proof of our algorithm, we use a result from \cite{Ram2012} which shows the averaged iterates can still arrive at consensus if the errors behave nicely.
\begin{lemma}\label{lemma:ram}
\textit{
Let Assumptions \ref{assume:ds} and \ref{assume:nc} hold.
%and such that $\sum_{k=0}^\infty \a_k=\infty$. NOT NEEDED
Consider the iterates generated by
\begin{equation}\label{eqn:rule}
\theta_{i,k} \hspace{-0.5mm}= \hspace{-0.5mm}\sum_{j\in\Vc} [W_k]_{ij} \theta_{j,k-1} + e_{i,k}, \text{ for } i \in \Vc.
\end{equation}
Let $\bar\theta_k$ denote the average of $\theta_{i,k}$ for $i \in \Vc$, i.e., $\bar\theta_k = \frac{1}{N}\sum_{i\in\Vc}\theta_{i,k}$.
Suppose there exists a nonnegative nonincreasing scalar sequence $\{\alpha_k\}$ such that
\[
\sum_{k=1}^{\infty} \alpha_k \|e_{i,k}\|< \infty, \text{ for all } i \in \Vc.
\]
Then, for all  $i,j \in \Vc$
\[
\sum_{k=1}^{\infty} \alpha_k \|\theta_{i,k}-\theta_{j,k}\|< \infty.
\]
Furthermore, for all $i\in \Vc$ and $k \ge 1$,
\begin{align*}
\|&\theta_{i,k}-\bar\theta_k\| \le N \gamma \beta^k \max_j \|\theta_{j,0}\|\\
&~+ \gamma\sum_{\ell=0}^{k-1}\beta^{k-\ell}\sum_{j=1}^N\|e_{j,\ell+1}\|+ \frac{1}{N}\sum_{j=1}^N\|e_{j,k}\|+\|e_{i,k}\|
\end{align*}
where $\gamma$ and $\beta$ are defined as
\begin{align*}
\gamma = \left(1-\frac{\nu}{2N^2}\right)^{-2}, \quad \beta = \left(1-\frac{\nu}{2N^2}\right)^{\frac{1}{Q}}.
\end{align*}
}
\end{lemma}
%Using this Lemma \ref{lemma:ram}, we estimate the disagreement among $p_i(k)$ for $i \in \Vc$ that will be used in our analysis later.

\subsection{Lemmas}
We need a series of lemmas for proving Proposition \ref{prop:as} and \ref{prop:as2}.
%The proofs of these lemmas are provided in the Appendix.
We first state an auxiliary lemma that will be later used to relate two consecutive iterates $x_{i,k}$ and $x_{i,k-1}$. This lemma can be shown by combining two existing results in \cite{Polyak-sg} and \cite{AN2011}, but we include it here for completeness.
\begin{lemma} \label{lem:xp}
\textit{
Let Assumptions \ref{assume:f} and \ref{assume:Pr} hold.
Let the iterates $\{p_{i,k}\}$, $\{v_{i,k}\}$ and $\{x_{i,k}\}$ be generated by the algorithm (\ref{eqn:algo1})-(\ref{eqn:algo3}).
Then, we have almost surely for any $\check{x},z \in \Xc_0$, $i \in \Vc$ and $k \ge 1$,
\begin{align*}
\|x_{i,k}-\check{x}\|^2 \le~& \|p_{i,k}-\check{x}\|^2 - 2\a_k(f_i(z)-f_i(\check{x}))\\
& - \frac{\tau-1}{\tau C_{g}^2} \left(g^+(p_{i,k},\o_{i,k})\right)^2\\
& + \frac{1}{4\eta}\|p_{i,k}-z\|^2 + D_{\tau,\eta}\a_k^2,
\end{align*}
where $D_{\tau,\eta} = (\tau+4\eta+1)C_{f_i}^2$ and $\eta, \tau >0$ are arbitrary.
}
\end{lemma}
\begin{IEEEproof}
In the light of \cite[Theorem 1]{Polyak-sg}, we obtain from algorithm (\ref{eqn:algo3}) and Assumption \ref{assume:f}(e)
\begin{equation}\label{eqn:polyak}
\|x_{i,k}-\check{x}\|^2 \le \|v_{i,k}-\check{x}\|^2 - \frac{\left(g^+(v_{i,k},\o_{i,k})\right)^2}{C_{g}^2},
\end{equation}
for any  $\check{x} \in \Xc_0$.
%From the uniform boundedness of the subgradients, we further obtain
%\[
%\|x-\check{x}\|^2 \le \|v-\check{x}\|^2 - \b(2-\b) \frac{\left(h^+(v)\right)^2}{C_h^2}.
%\]
We can rewrite $g^+(v_{i,k},\o_{i,k}) = \left(g^+(v_{i,k},\o_{i,k}) -g^+(p_{i,k},\o_{i,k})\right)+g^+(p_{i,k},\o_{i,k})$. Therefore,
\begin{align}\label{eqn:gfnpre}
&(g^+(v_{i,k},\o_{i,k}))^2 \nonumber\\
&\ge 2g^+(p_{i,k},\o_{i,k})\left(g^+(v_{i,k},\o_{i,k}) -g^+(p_{i,k},\o_{i,k})\right)\nonumber\\
&~~ + \left(g^+(p_{i,k},\o_{i,k})\right)^2.
\end{align}
The first term on the right-hand side of (\ref{eqn:gfnpre}) can be further estimated as
\begin{align}\label{eqn:gfn}
&2g^+(p_{i,k},\o_{i,k})\left(g^+(v_{i,k},\o_{i,k}) -g^+(p_{i,k},\o_{i,k})\right)\nonumber\\
&\ge-2g^+(p_{i,k},\o_{i,k})\left|g^+(v_{i,k},\o_{i,k}) -g^+(p_{i,k},\o_{i,k})\right| \nonumber\\
&\ge -2C_{g}\|v_{i,k}-p_{i,k}\|g^+(p_{i,k},\o_{i,k}),
\end{align}
where the last inequality is from relation (\ref{eqn:Cgi}).
From the definition in (\ref{eqn:algo2}) and Assumption \ref{assume:f}(c), we further have that
\begin{align}\label{eqn:gfnb}
2C_{g}&\|v_{i,k}-p_{i,k}\|g^+(p_{i,k},\o_{i,k}) \nonumber\\
&\le 2\a_k C_{g}C_{f_i}g^+(p_{i,k},\o_{i,k}) \nonumber\\
&\le \tau \a_k^2C_{g}^2C_{f_i}^2 + \frac{1}{\tau}\left(g^+(p_{i,k},\o_{i,k})\right)^2,
\end{align}
where the last inequality is obtained by using $2|a||b| \le \tau a^2 + \frac{1}{\tau} b^2$ and $\tau >0$ is arbitrary.
Using relations (\ref{eqn:gfn})-(\ref{eqn:gfnb}) in (\ref{eqn:gfnpre}), we obtain,
\begin{align*}
\big(g^+(&v_{i,k},\o_{i,k})\big)^2 \\
&\ge -\tau \a_k^2C_{g}^2C_{f_i}^2 + \left(1-\frac{1}{\tau}\right)\left(g^+(p_{i,k},\o_{i,k})\right)^2.
\end{align*}
Hence, for all $\check{x} \in \Xc_0$,
\begin{align}\label{eqn:partI}
\|x_{i,k}-\check{x}\|^2 \le &\|v_{i,k}-\check{x}\|^2 - \frac{\tau-1}{\tau C_{g}^2}\left(g^+(p_{i,k},\o_{i,k})\right)^2 \nonumber\\
&+\tau C_{f_i}^2 \a_k^2.
\end{align}
As the update rule in (\ref{eqn:algo2}) coincides with the algorithm in \cite{AN2011},
we can reuse another existing lemma \cite[Lemma 3]{AN2011}. That is, for any $\check{x},z\in \Xc_0$, we have
\begin{align}\label{eqn:nedich}
\|v_{i,k}-\check{x}& \|^2 \le \|p_{i,k}-\check{x}\|^2 - 2\a_k(f_i(z)-f_i(\check{x}))\nonumber\\
& + \frac{1}{4\eta}\|p_{i,k}-z\|^2 + \a_k^2(1+4\eta)C_{f_i}^2,
\end{align}
where $\eta >0 $ is arbitrary.
Substituting this inequality in relation (\ref{eqn:partI}) concludes the proof.
\end{IEEEproof}

Since we use an approximate projection, we cannot guarantee the feasibility of the iterates $\{x_{i,k}\}$ and $\{p_{i,k}\}$.  In the next lemma, we prove that $\{p_{i,k}\}$ and $\{x_{i,k}\}$ for all $i \in \Vc$ asymptotically achieve feasibility. To this end, we define the following quantity: For all $i\in \Vc$ and $k \ge 1$, $z_{i,k}$ is defined as the projection of $p_{i,k}$ on the feasible set $\Xc$, i.e.,
\begin{align}\label{eqn:zdef}
z_{i,k} \triangleq \mathsf{\Pi}_{\mathcal{X}}[p_{i,k}].
\end{align}
\begin{lemma}\label{lem:key1}
\textit{
Let Assumptions \ref{assume:ds} - \ref{assume:c} hold.
%Let the function $f_i$, for all $i \in \Vc$, be differentiable and $\nabla f_i(x)$ be Lipschitz continuous with constant $L$ over $\Rbb^n$.
Let the sequence $\{\a_k\}$ be such that $\sum_{k=1}^{\infty} \a_k^2 < \infty$.
Then, the iterates $\{p_{i,k}\}$ and $\{x_{i,k}\}$ generated by each agent $i \in \Vc$ via method (\ref{eqn:algo1})-(\ref{eqn:algo3}) satisfy:
\begin{itemize}
\item[(a)] $\displaystyle\sum_{k=1}^\infty \dist^2(p_{i,k},\Xc)<\infty \quad a.s.$
\item[(b)]
\smlr{$\displaystyle\sum_{k=1}^\infty \|x_{i,k}-z_{i,k}\|^2<\infty \quad a.s.$}
\end{itemize}
where $z_{i,k}$ is defined in \eqref{eqn:zdef}.
}
\end{lemma}
\begin{IEEEproof}
We use Lemma~\ref{lem:xp} with $\check{x}=z= z_{i,k}$.
Therefore, for any $\o_{i,k}\in \O_i$, $i\in \Vc$, and $k \ge 1$, we obtain
%\allowdisplaybreaks{
\begin{align}\label{eqn:lem2a}
 \|&x_{i,k}-z_{i,k}\|^2 \nonumber\\
\le~& \|p_{i,k}-z_{i,k}\|^2
 - \frac{\tau-1}{\tau C_g^2} \left(g^+(p_{i,k},\o_{i,k})\right)^2\nonumber\\
& + \frac{1}{4\eta}\|p_{i,k}-z_{i,k}\|^2 + D_{\tau,\eta}\a_k^2,
\end{align}
where $D_{\tau,\eta} = (\tau+4\eta+1)C_{f_i}^2$ and $\eta, \tau >0$ are arbitrary.
%where $A= (2+16c+2c\b(2-\b))L^2$ and $B = 2+8c+2c\b(2-\b)$.
By the definition of the projection, we have
\begin{align*}
\dist(p_{i,k},\Xc) =&~ \|p_{i,k} - z_{i,k}\|, \text{ and}\\
\dist(x_{i,k},\Xc)
=&~ \|x_{i,k}-\mathsf{\Pi}_{\mathcal{X}}[x_{i,k}]\| \le \|x_{i,k}-z_{i,k}\|.
\end{align*}
Upon substituting these estimates in relation (\ref{eqn:lem2a}), we obtain
\begin{align}\label{eqn:lem2b}
\dist^2(x_{i,k},\Xc) \le &~\dist^2(p_{i,k},\Xc)
- \frac{\tau -1}{\tau C_{g}^2} \left(g^+(p_{i,k},\o_{i,k})\right)^2\nonumber\\
&~+ \frac{1}{4\eta}\dist^2(p_{i,k},\Xc) + D_{\tau,\eta}\a_k^2.
\end{align}

%Let $\Fc_k$ denote the algorithm's history up to time $k$. i.e.,
%\[
%\Fc_k = \{x_{i,0}, ~(\o_{i,t},~1\le t \le k),~ i \in \Vc\},
%\]
%and $\Fc_= \{x_{i,0},~i\in \Vc\}$.
Taking the expectation conditioned on $\Fc_{k-1}$ and noting that $p_{i,k}$ is fully determined by $\Fc_{k-1}$,
we have almost surely for any $i \in \Vc$ and $k \ge 1$
\begin{align}\label{eqn:lem2c}
\Es\big[&\dist^2(x_{i,k},\Xc) \mid \Fc_{k-1}\big] \nonumber\\
\le~& \dist^2(p_{i,k},\Xc)
 - \frac{\tau-1}{\tau C_{g}^2} \Es\left[\left(g^+(p_{i,k},\o_{i,k})\right)^2 \mid \Fc_{k-1} \right]\nonumber\\
& + \frac{1}{4\eta}\dist^2(p_{i,k},\Xc) + D_{\tau,\eta}\a_k^2.
\end{align}
Furthermore, choosing $\tau =4$, $\eta = cC_{g}^2$ and using Assumption \ref{assume:c} yield
\begin{align}\label{eqn:lemkey1}
\Es&\left[\dist^2(x_{i,k},\Xc) \mid \Fc_{k-1}\right] \le \dist^2(p_{i,k},\Xc)\\
& - \frac{1}{2cC_{g}^2}\dist^2(p_{i,k},\Xc)
+ (5+4cC_{g}^2)C_{f_i}^2\a_k^2.\nonumber
\end{align}
Finally, by summing over all $i$ and using Lemma \ref{lem:ds} with $h(x) = \dist^2(x,\Xc)$,
we arrive at the following relation:
\begin{align} \label{eqn:lem2e}
\sum_{i \in \Vc}&\Es\left[\dist^2(x_{i,k},\Xc) \mid \Fc_{k-1}\right] \le \sum_{i \in \Vc} \dist^2(x_{i,k-1},\Xc) \nonumber\\
& - \frac{1}{2cC_{g}^2}\sum_{i \in \Vc}\dist^2(p_{i,k},\Xc)
+D N\a_k^2,
\end{align}
where $D = (5+4cC_g^2)C_f^2$ and $C_f = \max_{i\in\Vc} C_{f_i}$.
Therefore, for all $k \ge 1$, all the conditions of the convergence theorem (Theorem~\ref{thm:super}) are satisfied
and we conclude that
\begin{align}\label{eqn:psum}
\sum_{k=1}^\infty \dist^2(p_{i,k},\Xc)<\infty\quad\hbox{for all $i \in \Vc$ }\quad  a.s.
\end{align}
\smlr{Lastly, from relation \eqref{eqn:lem2a} and the chosen values for $\tau$ and $\eta$, we obtain
for any $i \in \Vc$ and $k \ge 1$ almost surely
\begin{align*}
\|x_{i,k}-z_{i,k}\|^2
\le \left(1+\frac{1}{4cC_g}\right)\dist^2(p_{i,k},\Xc) + D\a_k^2.
\end{align*}
Therefore, in view of the result in \eqref{eqn:psum} and $\sum_{k=1}^{\infty}\a_k^2 < \infty$, the relation above implies
\begin{align*}
\sum_{k=1}^\infty \|x_{i,k}-z_{i,k}\|^2<\infty\quad\hbox{for all $i \in \Vc$ }\quad  a.s.,
\end{align*}
which is our desired result.}
\end{IEEEproof}

To complete the proof, we show in part (a) of the next lemma that the error $e_{i,k}$
due to the perturbations made after the consensus step (\ref{eqn:algo1}), i.e.,
\begin{align}\label{eqn:edef}
e_{i,k}\triangleq x_{i,k}-p_{i,k},
\end{align}
eventually converges to zero for all $i \in \Vc$.
%This result combined with Lemma \ref{lem:key1} implies that the sequence $x_{i,k}$ for all $i\in \Vc$ also achieve asymptotic feasibility.
This will allow us to invoke Lemma \ref{lemma:ram} and show the iterate consensus.
In part (b) of the next lemma, we show that the sequences $\{z_{i,k}\}$ arrive at consensus by converging to their mean $\bar z_k$, i.e., for $k \ge 1$
\begin{align}\label{eqn:barzdef}
\bar z_k \triangleq \frac{1}{N}\sum_{i \in \Vc}z_{i,k}.
\end{align}
In part (c) of the next lemma, we show the network error term is summable.
%$\|z_{i,k}-\bar{z}_k\|$ for any $i\in \Vc$.
\begin{lemma} \label{lem:disagree}
\textit{
Let Assumptions \ref{assume:ds}-\ref{assume:nc} hold. \smlr{Let the sequence
$\{\a_k\}$ be nonnegative nonincreasing and $\sum_{k=1}^\infty \a_k^2<\infty$.}
%Define $\bar z_k \triangleq \frac{1}{N}\sum_{i \in \Vc}z_{i,k}$ and $e_{i,k} \triangleq x_{i,k}-p_{i,k}$ for all $i\in \Vc$ and $k\ge 1$.
Then, we have for all $i \in \Vc$
\begin{itemize}
\item[(a)]
$\displaystyle\sum_{k=1}^\infty \|e_{i,k}\|^2 <\infty \quad a.s.$
\item[(b)]
\smlr{$\lim_{k\to\infty}\|z_{i,k}-\bar z_k\| = 0\quad a.s.$}
\item[(c)]
$\displaystyle\sum_{k=1}^\infty\a_k\|z_{i,k}-\bar z_k\|<\infty    \quad a.s.$
\end{itemize}
}
\end{lemma}

\begin{IEEEproof}
\noindent {\it Part (a): }
From the relation (\ref{eqn:algo1})-(\ref{eqn:algo3}), $e_{i,k}$ in \eqref{eqn:edef} can be viewed as the perturbation that we make on $p_{i,k}$ after the network consensus step (\ref{eqn:algo1}).
\smlr{Consider $\|e_{i,k}\|$, for which we can write
\begin{align*}
\|e_{i,k}\|
\le{}& \|x_{i,k} - z_{i,k}\| +\|z_{i,k} - p_{i,k}\|.
\end{align*}
Applying $(a+b)^2 \le 2a^2 + 2b^2$ in the above inequality, we have
\begin{align*}
\|e_{i,k}\|^2 \leq 2\|x_{i,k} - z_{i,k}\|^2 +2\dist^2(p_{i,k},\Xc).
\end{align*}
Summing this over $k$ and using Lemma \ref{lem:key1}, we obtain the desired result. }

\smlr{\noindent {\it Part (b): }
By applying the inequality $2ab\le a^2+b^2$ to each term in $\a_k\|e_{i,k}\|$ and using Lemma \ref{lem:disagree}(a), we further obtain for all $i \in \Vc$
\begin{align}\label{eqn:lem_e}
\sum_{k=1}^\infty \a_k \|e_{i,k}\|\le
\frac{1}{2}\sum_{k=1}^\infty\a_k^2 + {}&\frac{1}{2}\sum_{k=1}^\infty \|e_{i,k}\|^2<\infty
\quad a.s.
\end{align}
Using the relation above, \eqref{eqn:algo1} and $x_{i,k}=p_{i,k}+e_{i,k}$,
we can invoke Lemma~\ref{lemma:ram} with
$\theta_{i,k} = x_{i,k}$.
Therefore, it follows that
\begin{equation}\label{eqn:sumfin}
\sum_{k=1}^\infty\a_k\|x_{i,k}-x_{j,k}\| <\infty \hbox{ for all } i,j \in \Vc \quad a.s.
\end{equation}
Furthermore, for all $i\in \Vc$ and $k \ge 1$,
\begin{align}\label{eqn:xrec}
\|&x_{i,k}-\bar x_k\| \le N \gamma \beta^k \max_j \|x_{j,0}\|\\
&~+ \gamma\sum_{\ell=0}^{k-1}\beta^{k-\ell}\sum_{j=1}^N\|e_{j,\ell+1}\|+ \frac{1}{N}\sum_{j=1}^N\|e_{j,k}\|+\|e_{i,k}\|.\nonumber
\end{align}
From the fact that $0<\beta<1$ and part(a), we know the first term and the last two terms on the right-hand side converge to zero.
To show the second term also converges to zero,
we will use the following result from \cite[Lemma 3.1(a)]{Ram2010}.
\begin{lemma}
Let $\zeta_k$ be a scalar sequence. If $\lim_{k\to\infty} \zeta_k = \zeta$ and $0<\beta<1$,
then $\lim_{k\to\infty} \sum_{\ell=0}^{k-1}\beta^{k-\ell}\zeta_{\ell} = \frac{\zeta}{1-\beta}$.
\end{lemma}
From this lemma and the result in part(a), we know that the second term on the right-hand side of \eqref{eqn:xrec} also converges to zero. Therefore, we have for all $i \in \Vc$
\begin{align}\label{eqn:xcon}
\lim_{k\to\infty}\|x_{i,k} - \bar x_k\| = 0.
\end{align}
 We next consider the term $\|z_{i,k}-\bar{z}_k\|$, for which by using
 $\bar{z}_k = \frac{1}{N}\sum_{\ell \in\Vc} z_{\ell,k}$ we have
\begin{align*}
\|z_{i,k}-\bar{z}_k&\|=\left\|\frac{1}{N}\sum_{\ell\in\Vc}(z_{i,k}-z_{\ell,k})\right\|\cr
\le{}& \frac{1}{N}\sum_{\ell\in\Vc}\|z_{i,k}-z_{\ell,k}\|
 \le \frac{1}{N}\sum_{\ell\in\Vc}\|p_{i,k}-p_{\ell,k}\|,\nonumber
\end{align*}
 where the first inequality is obtained by the convexity of the norm
 and the last inequality follows by the non-expansive projection property in Lemma \ref{lem:proj}.
Furthermore, by using $\|p_{i,k}-p_{\ell,k}\|\le \|p_{i,k}-\bar p_k\|+ \|p_{\ell,k}- \bar p_k\|$, we obtain
for every $i\in \Vc$
\begin{align}\label{eqn:zp}
\|z_{i,k}-\bar{z}_k\|\le \|p_{i,k}-\bar p_k\|+\frac{1}{N}\sum_{\ell\in\Vc}\|p_{\ell,k}-\bar p_k\|.
\end{align}
We next consider $\|p_{i,k}- \bar p_k\|$.
By using the doubly stochasticity of $W_k$,
convexity of the norm and the fact that $0 \le [W_k]_{ij}\le 1$, we obtain
\begin{align*}
\|p_{i,k}- \bar p_k\|\le &\sum_{j\in\Vc} [W_k]_{ij}\left\|x_{j,k-1}-\bar p_k\right\|\\
\le & \sum_{j\in\Vc} \left\|x_{j,k-1}- \frac{1}{N}\sum_{\ell\in\Vc} x_{\ell,k-1}\right\|,
\end{align*}
where in the last equality we use
$\bar p_k= \frac{1}{N}\sum_{\ell \in \Vc}\sum_{j\in\Vc} [W_k]_{\ell j} x_{j,k-1} = \frac{1}{N}\sum_{\ell\in\Vc} x_{\ell,k-1}$.
Therefore, by using the convexity of the norm again, we see
\begin{align}\label{eqn:pmidx}
&\|p_{i,k}- \bar p_k\| \le \frac{1}{N}\sum_{j\in\Vc} \sum_{\ell\in\Vc}\left\|x_{j,k-1}- x_{\ell,k-1}\right\|\\
&~~~\le \frac{1}{N}\sum_{j\in\Vc} \sum_{\ell\in\Vc}\left(\|x_{j,k-1}- \bar x_{k-1}\|+\|x_{\ell,k-1}- \bar x_{k-1}\|\right).\nonumber
\end{align}
Combining this relation with \eqref{eqn:zp} and using the result in \eqref{eqn:xcon}, we obtain the desired result.
}

\noindent {\it Part (c): }
By using relation~\eqref{eqn:sumfin} in \eqref{eqn:pmidx}, we obtain
\begin{align}\label{eqn:lessinf}
&\sum_{k=1}^{\infty}\a_k\|p_{i,k}- \bar p_k\|\nonumber\\
&\le \sum_{k=1}^{\infty}\frac{\a_k}{N}\sum_{j\in\Vc} \sum_{\ell\in\Vc}\left\|x_{j,k-1}- x_{\ell,k-1}\right\| < \infty.
\end{align}
Upon summing the relation \eqref{eqn:zp} over $i\in \Vc$, we find
\begin{align}\label{eqn:sums}
\sum_{i \in \Vc} \|z_{i,k}-\bar{z}_k\|\le 2 \sum_{i\in\Vc} \|p_{i,k}-\bar p_k\|.
\end{align}
Therefore, from (\ref{eqn:lessinf}) and (\ref{eqn:sums}), we obtain
\begin{align*}
\sum_{i \in \Vc} \sum_{k=1}^{\infty} \a_k\|z_{i,k}-\bar{z}_k\|\le 2  \sum_{i\in\Vc} \sum_{k=1}^{\infty}\a_k\|p_{i,k}-\bar p_k\| < \infty,
\end{align*}
which is the desired result.
\end{IEEEproof}

In the next lemma, we use standard convexity analysis to lower-bound the term $\sum_{i \in \Vc}(f_i(z_{i,k})-f_i(\check{x}))$ with a network error term and a global term.
\begin{lemma}\label{lem:disagree_f}
\textit{
Let Assumption \ref{assume:f} hold. Then, for all $\check{x} \in \Xc$, we have
\[
\sum_{i \in \Vc}\hspace{-0.7mm}(f_i(z_{i,k})-f_i(\check{x})) \hspace{-0.7mm}\ge \hspace{-0.7mm}-C_f \sum_{i \in \Vc}\hspace{-0.7mm} \|z_{i,k}-\bar{z}_k\| + f(\bar{z}_k) - f(\check{x}),
\]
where $C_f = \max_{i\in\Vc} C_{f_i}$.
}
\end{lemma}

\begin{IEEEproof}
Recall that $f(x)=\sum_{i=1}^m f_i(x)$. Recall that $\bar{z}_k = \frac{1}{N}\sum_{\ell \in \Vc} z_{\ell,k}$.
Using $\bar{z}_k$ and $f$,
we can rewrite the term $f_i(z_{i,k})-f_i(\bar{x})$ as follows:
\begin{align}\label{eqn:rewrite}
\sum_{i\in\Vc} &(f_i(z_{i,k})-f_i(\check{x})) \nonumber\\
 &=  \sum_{i\in\Vc} (f_i(z_{i,k})-f_i(\bar{z}_k)) +  (f(\bar{z}_k)-f(\check{x})).
\end{align}
Furthermore, using the convexity of each function $f_i$, we obtain
\begin{align*}
\sum_{i\in \Vc} (f_i(z_{i,k})-f_i(\bar{z}_k))
& \geq \sum_{i\in\Vc}  \langle s_{i,k},z_{i,k}-\bar{z}_k\rangle\\
& \ge  -\sum_{i\in\Vc}  \|s_{i,k}\|\,\|z_{i,k}-\bar{z}_k\|,
\end{align*}
where $s_{i,k}$ is a subgradient of $f_i$ at $\bar{z}_k$.
Since $\bar{z}_k$ is a convex combination of points $z_{i,k}\in \Xc \subseteq \Xc_0$, it follows that
$\bar z_k \in \Xc_0$.
This observation and Assumption~\ref{assume:f}(c), stating that the subgradients of $f_i(x)$ are uniformly bounded for $x\in\Xc_0$, yield
\begin{align}\label{eqn:a}
\sum_{i\in\Vc} (f_i(z_{i,k})-f_i(\bar{z}_k))
\ge -C_f\sum_{i\in\Vc} \|z_{i,k}-\bar{z}_k\|,
\end{align}
where $C_f = \max_{i\in\Vc}C_{f_i}$.
Therefore, from (\ref{eqn:rewrite}) and (\ref{eqn:a}), we have that
\begin{align*}
\sum_{i \in \Vc}&(f_i(z_{i,k})-f_i(\check{x})) \\
&\ge -C_f \sum_{i \in \Vc} \|z_{i,k}-\bar{z}_k\| + f(\bar{z}_k) - f(\check{x}).
\end{align*}
\end{IEEEproof}

\subsection{Proof of Proposition \ref{prop:as}}
We invoke Lemma \ref{lem:xp} with $z = z_{i,k} = \proj_{\Xc}[p_{i,k}]$, $\tau = 4$ and $\eta = cC_{g}^2$. We also let $\check{x} = x^*$ for an arbitrary $x^* \in \Xc^*$. Therefore, for any $x^* \in \Xc^*$, $i \in \Vc$ and $k \ge 1$, we almost surely have
\begin{align}\label{eqn:propmid}
\|x_{i,k}-&x^*\|^2 \le \|p_{i,k}-x^*\|^2 - 2\a_k(f_i(z_{i,k})-f_i(x^*)) \\
& - \frac{3}{4C_{g}^2} \left(g^+(p_{i,k},\o_{i,k})\right)^2 \nonumber\\
&+ \frac{1}{4cC_{g}^2}\dist^2(p_{i,k},\Xc)
 + (5+4cC_{g}^2)C_{f_i}^2\a_k^2. \nonumber
\end{align}
%Let $\Fc_k$ denote the algorithm's history up to time $k$. i.e.,
%\[
%\Fc_k = \{x_{i,0}, ~(\o_{i,t},~1\le t \le k),~ i \in \Vc\},
%\]
%and $\Fc_= \{x_{i,0},~i\in \Vc\}$.
By taking the expectation conditioned on $\Fc_{k-1}$ in the above relation and summing this over $i \in \Vc$, we obtain \begin{align*}
\sum_{i\in \Vc}&\Es\left[\|x_{i,k}-x^*\|^2 \mid \Fc_{k-1}\right] \\
~&\le \sum_{i\in \Vc}\|p_{i,k}-x^*\|^2 - 2\a_k\sum_{i\in \Vc}(f_i(z_{i,k})-f_i(x^*)) \\
~& - \frac{3}{4 C_g^2} \sum_{i\in \Vc}\Es\left[\left(g^+(p_{i,k},\o_{i,k})\right)^2 \mid \Fc_{k-1}\right]\\
~&+ \frac{1}{4cC_g^2}\sum_{i\in \Vc}\dist^2(p_{i,k},\Xc)
 + DN\a_k^2,
\end{align*}
where $D = (5+4cC_g^2)C_f^2$ with $C_f = \max_{i\in \Vc} C_{f_i}$.
Now we use Lemma \ref{lem:ds} with $h(x) = \|x-x^*\|^2$, Assumption \ref{assume:c} and Lemma \ref{lem:disagree_f} with $\check{x} = x^*$ to further estimate the terms on the right-hand side. From these, obtain almost surely for any  $k \ge 1$ and $x^* \in \Xc^*$,
\begin{align*}
\sum_{i \in \Vc}&\Es\left[\|x_{i,k}-x^*\|^2 \mid \Fc_{k-1}\right] \\
&\le \sum_{i \in \Vc}\|x_{i,k-1}-x^*\|^2 - 2\a_k\sum_{i \in \Vc}(f(\bar{z}_k)-f(x^*)) \\
& - \frac{1}{2cC_g^2} \sum_{i \in \Vc}\dist^2(p_{i,k},\Xc) \nonumber\\
& + 2\a_k C_f \sum_{i \in \Vc} \|z_{i,k}-\bar{z}_k\| + DN\a_k^2. \nonumber
\end{align*}
Since $\bar{z}_k \in \Xc$, we have $f(\bar{z}_k) - f(x^*) \ge 0$.
Thus, under the assumption $\sum_{k=0}^{\infty} \alpha_k^2 < \infty$ and Lemma \ref{lem:disagree}(c), the above relation satisfies all the conditions of the convergence
Theorem~\ref{thm:super}.
\smlr{
Using this theorem, we have the following results.}

\smlr{
\textit{Result 1:}  The sequence $\{\sum_{i\in\Vc}\|x_{i,k}-x^*\|\}$ is
convergent \textit{a.s.} for every $x^* \in \Xc^*$.}

\smlr{
\textit{Result 2:}   For every $x^* \in \Xc^*$,
\[
\sum_{k=1}^{\infty} \alpha_k(f(\bar{z}_k)-f(x^*)) < \infty \quad a.s.
\]
}

\smlr{
From \textit{Result 1} and Lemma \ref{lem:key1}(b),
we know that the sequence $\{\sum_{i\in\Vc}\|z_{i,k}-x^*\|\}$ is
convergent \textit{a.s.} for every $x^* \in \Xc^*$.
This and Lemma \ref{lem:disagree}(b) imply that $\|\bar z_k - x^*\|$
is also convergent \textit{a.s.} for every $x^* \in \Xc^*$.
From \textit{Result 2}, $\sum_{k=1}^{\infty}\a_k = \infty$, and the continuity of $f$, it follows that
the sequence $\{\bar{z}_k\}$ must have one accumulation point in the set $\Xc^*$ \textit{a.s.}
This and the fact that $\{\|\bar{z}_k-x^*\|\}$ is convergent \textit{a.s.} for every $x^* \in \Xc^*$ imply that for a random point $x^\star \in \Xc^*$,
\begin{equation}\label{eqn:z_final1}
\lim_{k \rightarrow \infty} \bar{z}_k = x^\star \quad a.s.
\end{equation}
}

\smlr{
We now prove the following claim: For all $i \in \Vc$
\begin{equation}\label{eqn:x_final1}
\lim_{k \rightarrow \infty} x_{i,k} = x^\star \quad a.s.
\end{equation}
Consider
\begin{align*}
\|x_{i,k}-x^\star\| \le \|x_{i,k}-z_{i,k}\| + \|z_{i,k}-\bar z_k\| + \|\bar z_k -x^\star\|.
\end{align*}
From Lemma \ref{lem:key1}(b), Lemma \ref{lem:disagree}(b) and \eqref{eqn:z_final1},
all the terms on the right-hand side converge to zero \textit{a.s.}
Therefore, it is obvious that claim \eqref{eqn:x_final1} holds, which is our desired result.
}

\subsection{Proof of Proposition \ref{prop:as2}}
The line of proof is similar to that in Proposition \ref{prop:as}. Therefore, we only lay down the differences.

Using the definition of $\pi$ in Lemma \ref{lem:perron}, we have
\begin{align}\label{eqn:dsToP}
[\pi]_j = \sum_{i\in\Vc}[\pi]_i[W]_{ij}, \text{ for all } j \in \Vc.
\end{align}
Also, in the proof we consider the following weighted averages rather than the true averages $\bar x_k$, $\bar p_k$ and $\bar z_k$.
\begin{align}
&\hat x_k \triangleq \sum_{i\in\Vc} [\pi]_i x_{i,k},\quad
\hat p_k \triangleq \sum_{i\in\Vc} [\pi]_i p_{i,k},\\
&\text{ and } \hat z_k \triangleq \sum_{i\in\Vc} [\pi]_i z_{i,k}.\nonumber
\end{align}

First, notice that Lemma \ref{lem:xp} still holds in this case as it does not require Assumption \ref{assume:ds}.\\

\noindent\textit{Changes in Lemma \ref{lemma:ram}: }
Combining with the results in \cite{Blondel05convergencein,con01,Moreau_Stability,ANAO2010,Tsi1986,ANAO}, Lemma \ref{lemma:ram} still holds in this case by replacing $\bar \theta_k$ with $\hat \theta_k \triangleq \sum_{i\in\Vc} [\pi]_i \theta_{i,k}$ and re-defining the constants $\gamma$ and $\beta$ as
\begin{align*}
\gamma = 2, \quad \beta = 1-\frac{1}{N^N}.
\end{align*}
If in addition the underlying graph $\Gc$ is regular, then we have
\begin{align*}
\gamma = \sqrt{2}, \quad \beta = \min\left\{1-\frac{1}{4N^3},~\max_{k\ge1} \sigma_2(W)\right\},
\end{align*}
where $\sigma_2(W)$ is the second largest singular value of $W$.
\\

\noindent\textit{Changes in Lemma \ref{lem:key1}: }
By multiplying $[\pi]_i$ to \eqref{eqn:lemkey1} and summing over $i \in \Vc$, we obtain
\begin{align} \label{eqn:keydiff}
\sum_{i \in \Vc}&\Es\left[[\pi]_i\dist^2(x_{i,k},\Xc) \mid \Fc_{k-1}\right] \le \sum_{i \in \Vc} [\pi]_i\dist^2(p_{i,k},\Xc) \nonumber\\
& - \frac{1}{2cC_{g}^2}\sum_{i \in \Vc}[\pi]_i\dist^2(p_{i,k},\Xc)
+D N\a_k^2.
\end{align}
From the definition of $p_{i,k}$ in \eqref{eqn:algo1} and the convexity of the distance function, we have
\begin{align*}
\sum_{i \in \Vc} [\pi]_i\dist^2(p_{i,k},\Xc) \le &~ \sum_{i \in \Vc} \sum_{j\in \Vc}[\pi]_i[W]_{ij}\dist^2(x_{j,k-1},\Xc)\\
\le &~\sum_{j\in \Vc}[\pi]_j\dist^2(x_{j,k-1},\Xc),
\end{align*}
where the last inequality follows from \eqref{eqn:dsToP}.
Combining this result with \eqref{eqn:keydiff}, we obtain
\begin{align*}
\sum_{i \in \Vc}&\Es\left[[\pi]_i\dist^2(x_{i,k},\Xc) \mid \Fc_{k-1}\right] \le \sum_{i \in \Vc} [\pi]_i\dist^2(x_{i,k-1},\Xc) \nonumber\\
& - \frac{1}{2cC_{g}^2}\sum_{i \in \Vc}[\pi]_i\dist^2(p_{i,k},\Xc)
+D N\a_k^2,
\end{align*}
in which all the conditions of Theorem \ref{thm:super} holds. Hence, all the remaining results follow immediately.\\

\noindent\textit{Changes in Lemma \ref{lem:disagree}: }
All the results still hold by replacing $\bar z_k$, $\bar p_k$ and $\bar x_k$ with $\hat z_k$, $\hat p_k$ and $\hat x_k$, respectively. Especially, from relation \eqref{eqn:dsToP} we have
\begin{align*}
\hat p_k = \sum_{i\in\Vc} [\pi]_i \sum_{j\in \Vc} [W]_{ij}x_{j,k-1} = \sum_{j\in \Vc} [\pi]_j x_{j,k-1} = \hat x_{k-1},
\end{align*}
and all the results follow immediately.
\\

\noindent\textit{Changes in Proposition \ref{prop:as}: }
By multiplying $[\pi]_i$ to \eqref{eqn:propmid}, summing this over $i \in \Vc$
and considering $f_i(x) = \frac{1}{N} a^{\top}x$, we have
\begin{align*}
&\sum_{i\in\Vc}[\pi]_i\|x_{i,k}-x^*\|^2 \le \sum_{i\in\Vc}[\pi]_i\|p_{i,k}-x^*\|^2 \\
&~~~- \frac{2\a_k}{N}a^{\top}(\hat z_k-x^*)
 - \frac{3}{4C_{g}^2} \sum_{i\in\Vc}[\pi]_i\left(g^+(p_{i,k},\o_{i,k})\right)^2 \nonumber\\
&~~~+ \frac{1}{4cC_{g}^2}\sum_{i\in\Vc}[\pi]_i\dist^2(p_{i,k},\Xc)
 + DN\a_k^2, \nonumber
\end{align*}
where we used the fact that $\sum_{i\in\Vc} [\pi]_i z_{i,k} = \hat z_k$.

Now we use \eqref{eqn:dsToP} and Assumption \ref{assume:c} to obtain almost surely for any  $k \ge 1$ and $x^* \in \Xc^*$,
\begin{align*}
\sum_{i \in \Vc}&\Es\left[[\pi]_i\|x_{i,k}-x^*\|^2 \mid \Fc_{k-1}\right] \\
&\le \sum_{i \in \Vc}[\pi]_i\|x_{i,k-1}-x^*\|^2 - \frac{2\a_k}{N} a^{\top}(\hat z_k-x^*) \\
& - \frac{1}{2cC_g^2} \sum_{i \in \Vc}[\pi]_i\dist^2(p_{i,k},\Xc)
 + DN\a_k^2. \nonumber
\end{align*}
Since $\hat{z}_k \in \Xc$, we have $a^{\top}(\hat z_k-x^*) \ge 0$.
Thus, under the assumption $\sum_{k=0}^{\infty} \alpha_k^2 < \infty$ and Lemma \ref{lem:key1}(a), the above relation satisfies all the conditions of the convergence Theorem~\ref{thm:super}.

Using this theorem, we have the following results.

\textit{Result 1:}  The sequence $\{\sum_{i\in\Vc}[\pi]_i\|x_{i,k}-x^*\|\}$ is
convergent \textit{a.s.} for every $x^* \in \Xc^*$.

\textit{Result 2:}   For every $x^* \in \Xc^*$,
\[
\sum_{k=1}^{\infty} \alpha_ka^{\top}(\hat z_k-x^*) < \infty \quad a.s.
\]

From \textit{Result 1} and Lemma \ref{lem:key1}(b),
we know that the sequence $\{\sum_{i\in\Vc}[\pi]_i\|z_{i,k}-x^*\|\}$ is
convergent \textit{a.s.} for every $x^* \in \Xc^*$.
This and Lemma \ref{lem:disagree}(b) imply that $\|\hat z_k - x^*\|$
is also convergent \textit{a.s.} for every $x^* \in \Xc^*$.
From \textit{Result 2}, $\sum_{k=1}^{\infty}\a_k = \infty$, it follows that
the sequence $\{\hat{z}_k\}$ must have one accumulation point in the set $\Xc^*$ \textit{a.s.}
This and the fact that $\{\|\hat{z}_k-x^*\|\}$ is convergent \textit{a.s.} for every $x^* \in \Xc^*$ imply that for a random point $x^\star \in \Xc^*$,
\begin{align*}
\lim_{k \rightarrow \infty} \hat{z}_k = x^\star \quad a.s.
\end{align*}
The remaining results follow immediately.

\section{Simulation Results\label{sec:sim}}
In this section, we provide a numerical example showing the effectiveness of the proposed decentralized approximate projection algorithm.
We consider optimal gossip averaging which is an example of decentralized optimization.
%In Section \ref{sec:lqr},
%we consider robust LQR which is a feasibility problem that can be solved in a decentralized fashion.

%\subsection{Optimal Gossip Averaging \label{sec:gossip}}
In many decentralized algorithms, gossip based communication protocols are often used.
In these communication protocols,
only one agent randomly wakes up at a time (say agent $i$) and selects one of its neighbors (say agent $j$) with probability $p_{ij}$. Then, the two agents exchange their current information through the link $(i,j)$ and take the average.
Let $A(i,j)$ denote the averaging matrix associated with the link $(i,j)$. For example, the averaging matrix $A(1,2)$ of a 4-agent network system looks like
\[
A(1,2) = \left[
\begin{array}{cccc}
1/2& 1/2 & 0 & 0 \\
1/2& 1/2 & 0 & 0 \\
0& 0 & 1 & 0 \\
0& 0 & 0 & 1
\end{array}
\right].
\]
Note that the expectation of the averaging matrix $A$ can be represented as $\Es[A]= \frac{1}{N}\sum_{i,j\in\Vc}p_{ij}A(i,j)$.

Let $P$ denote the probability matrix whose component of the $i$-th row and $j$-th column is $p_{ij}$.
%, which is the probability of agent $i$ choosing agent $j$ among its neighbors.
Our goal here is to find an optimal probability matrix $P^*$
associated with the current communication graph, which is time-invariant and connected,  in a decentralized fashion.
%As this is a gossip averaging problem, the graph associated with $P$ must be undirected and assumed to be connected.
The convergence speed of the gossip protocol is known to be inversely proportional to $\lambda_2(\Es(A))$, which is the second largest eigenvalue of the expected averaging matrix $\Es[A]$ (see \cite{Boyd:2006}).
Thus, the optimization problem of finding the fastest averaging distribution $P^*$ can be formulated as the following SDP:
\begin{subequations}
\begin{align}
\min_{s,P}~ & s \label{prob:gossip1}\\
\text{s.t. } & \sum_{i,j\in\Vc}p_{ij}A(i,j) - \1b\1b^{\top} \preceq sI \label{prob:gossip2}\\
& p_{ij} \ge 0, \quad p_{ij} = 0 \text{ if } (i,j) \not\in \Ec\label{prob:gossip3}\\
& \sum_{j\in\Vc} p_{ij} = 1, \text{ for } i \in \Vc. \label{prob:gossip4}
\end{align}
\end{subequations}
An optimal $P^*$ of the problem (\ref{prob:gossip1})-(\ref{prob:gossip4}) computed in a centralized fashion is not useful as gossip protocol is usually required in a decentralized setting.
A decentralized method has been proposed in \cite{Boyd:2006}, but the method only finds a suboptimal solution.
%of a reformulated problem.
%This is because the SDP constraint (\ref{prob:gossip2}) is not easily decentralized.
Using our proposed algorithm, we can find the optimal solution of (\ref{prob:gossip1})-(\ref{prob:gossip4}) in a decentralized way.

With a slight abuse of notation, let $x \triangleq [s; P]$. In this problem, all agents share the same local objective function, i.e., $f_i(x) = s$ for all $i\in\Vc$, whereas each agent $i$ has a local constraint set $\Xc_i = \Xc_i^1 \medcap \Xc_i^2$
where
\[
\Xc_i^1 = \{x\mid \sum_{i,j\in\Vc}p_{ij}A(i,j) - \1b\1b^{\top} \preceq sI\},
\]
\[
\Xc_i^2 = \{x\mid \sum_{j\in\Vc} p_{ij} = 1,~
p_{ij} \ge 0,~ p_{ij} = 0 \text{ if } j \not\in \Nc_i \}.
\]
At each iteration of our algorithm, we randomly select a component constraint from $\Xc_i$ and make a projection.
More specifically, we approximate the projection onto the SDP constraint $\Xc_i^1$ using the equation (\ref{eqn:subgradLMI}).
Note that the constraints \eqref{prob:gossip3}-\eqref{prob:gossip4}, which are distributed among agents, will guarantee the structure of the underlying communication graph. Therefore, agents do not require knowledge on the whole graph structure.
%The matrix $P$ here is just an optimization variable, where each agent makes a local copy.

\sml{We note that due to the compactness of the set $\Xc_i^2$, the problem \eqref{prob:gossip1}-\eqref{prob:gossip4} satisfies Assumption \ref{assume:f} and the optimal solution set $\Xc^*$ is nonempty.
Assumptions \ref{assume:ds}, \ref{assume:Pr} and \ref{assume:nc} can be satisfied by construction.
Assumption \ref{assume:c} is also satisfied as all inequalities are affine in this case.
%usually not easy to show in advance, but its sufficient conditions (e.g., nonempty interior of the feasible set or affine inequalities) are general enough to cover many optimization problems, so we assume the problem (\ref{prob:gossip1})-(\ref{prob:gossip4}) satisfy at least one of these sufficient conditions.
}

We let all agents terminate if their solution is within 0.01\% of the global average
and the total feasibility violation is less than $0.001$.
We say the algorithm has converged only when all network agents terminate.
Note that this global average based criterion is just used for the sake of simulations.
%In a real setting, we can change the termination criterion in such a way that each agent stops based on its local average.For example, each agent can keep track of the iterates of its neighboring agents and terminate if its own solution is within 0.01\% of the local average.
\sml{Also, due to the randomness of our algorithm, we repeat all the simulations for 10 times and report their averages.}

Table \ref{tbl:gossip} summarizes the simulation results. It shows the number of iterations until convergence for different numbers of agents ($N$) and underlying communication topologies ($\Gc_k$).
In the experiment, we use $\a_{i,k} = \frac{1}{k}$ for all $i \in \Vc$, and $4, 15$ agents with three different network topologies, namely clique, cycle and star.
\sml{Note that for this problem the underlying network must be time-invariant, i.e., $\Gc_k = \Gc$ for $k \ge 1$,
as the gossip algorithm in \cite{Boyd:2006} is built on a fixed undirected graph.}
%For simplicity, we assume the underlying network is time-invariant, i.e., $\Gc_k = \Gc$ for $k \ge 1$.
As expected, the star graph takes the most iterations for both $N = 4, 15$.
Also, when there are more agents in the network, the algorithm takes more iterations.

\begin{comment}
0.1%
	clique		cycle 		line
4	0.1660 (221)	0.2491 (287)	0.2805 (944)	
20	0.0470 (222)	0.0971(1026)	0.1019  (7096)

0.01%
	clique		cycle 		line		star
4	0.1666 (2170)	0.2499 (2819)	0.2807 (7607)	0.2742 (7190)
15	0.0619 (2179)	0.1275 (7018)	0.1354 (>50000)	0.4431 (18541)
20	0.0473 (1,984)	0.0975 (8280)	0.1020 (>50000)
\end{comment}

\begin{center}
\begin{table}
\caption{Number of iterations for all network agents to converge within 0.01\% of the global average \label{tbl:gossip}}
\centering{
\begin{tabular}{|c||c|c|c|}\hline
& clique & cycle & star \\\hline
$N=4$ & \hfill 2,170& \hfill 2,819& \hfill 7,190\\\hline
$N=15$ &\hfill 2,179  &\hfill 8,280  &\hfill 18,541 \\\hline
\end{tabular}
}
\end{table}
\end{center}

%For a stopping criterion, we first obtain an optimal solution $x^*$ to the problem (\ref{prob:gossip1})-(\ref{prob:gossip4}) by solving it in a centralized fashion. The decentralized algorithm terminates if the iterates of all network agents converge within 0.1\% of $x^*$.

\section{Conclusion\label{sec:con}}
We have studied a distributed optimization problem defined on a multiagent network which involves nontrivial constraints like LMIs. We have proposed a decentralized algorithm based on random feasibility updates, where we approximate the projection with an additional subgradient step.
The proposed algorithm is efficiently applicable for solving any distributed optimization problems
which involve lots of computationally prohibitive constraints, for example, decentralized SDPs.
We have established the almost sure convergence of our method
under two different assumptions on the weight matrices,
namely doubly stochastic $\{W_k\}$ over a $Q$-strongly connected sequence of digraphs
and row stochastic $W$ over a strongly connected digraph.
We have performed experiments on an optimal gossip averaging problem to verify the performance and convergence of the proposed algorithm.

\bibliography{soomin002}{}
\bibliographystyle{IEEEtran}
%\bibliographystyle{plain}

%\begin{comment}
\begin{IEEEbiography}[{\includegraphics[width=1in,height=1.25in,clip,keepaspectratio]{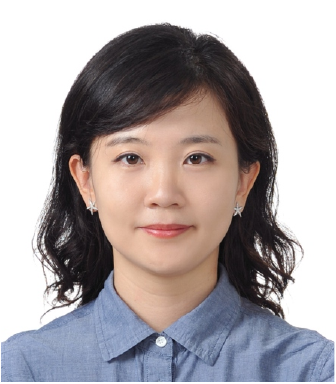}}]
{Soomin Lee} is currently working as a Research Scientist at Yahoo! Labs. She received her Ph.D. in Electrical and Computer Engineering from the University of Illinois, Urbana-Champaign (2013). She received two master's degrees from the Korea Advanced Institute of Science and Technology in Electrical Engineering, and from the University of Illinois at Urbana-Champaign in Computer Science.
After graduation, she worked as Postdoctoral Associate in Mechanical Engineering and Materials Science at Duke University and Postdoctoral Researcher in Industrial and Systems Engineering at Georgia Tech.  
%In 2009, she was an assistant research officer at the Advanced Digital Science Center (ADSC) in Singapore. 
Her research interests include theoretical optimization (convex, non-convex, online and stochastic), distributed control and optimization of various engineering systems interconnected over complex networks and large-scale machine learning in Internet industry.
\end{IEEEbiography}

\begin{IEEEbiography}[{\includegraphics[width=1in,height=1.25in,clip,keepaspectratio]{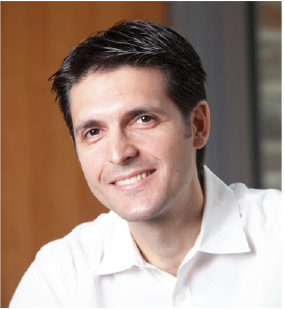}}]
{Michael M. Zavlanos} (S’05–M’09) received the Diploma in mechanical engineering from the National Technical University of Athens (NTUA), Athens, Greece, in 2002, and the M.S.E. and Ph.D. degrees in electrical and systems engineering from the University of Pennsylvania, Philadelphia, PA, in 2005 and 2008, respectively.

He is currently an Assistant Professor in the Department of Mechanical Engineering and Materials Science at Duke University, Durham, NC. He also holds a secondary appointment in the Department of Electrical and Computer Engineering and the Department of Computer Science. Prior to joining Duke University, Dr. Zavlanos was an Assistant Professor in the Department of Mechanical Engineering at Stevens Institute of Technology, Hoboken, NJ, and a Postdoctoral Researcher in the GRASP Lab, University of Pennsylvania, Philadelphia, PA. His research interests include a wide range of topics in the emerging discipline of networked systems, with applications in robotic, sensor, and communication networks. He is particularly interested in hybrid solution techniques, on the interface of control theory, distributed optimization, estimation, and networking.

Dr. Zavlanos is a recipient of various awards including the 2014 Office of Naval Research Young Investigator Program (YIP) Award and the 2011 National Science Foundation Faculty Early Career Development (CAREER) Award.
\end{IEEEbiography}
%\end{comment}

\end{document}